\documentclass[reqno,showkeys]{amsart}

\usepackage{amsfonts}
\usepackage{amssymb}
\usepackage{amsthm}
\usepackage{upref}

\makeatletter
 
\def\LaTeX{\leavevmode L\raise.42ex
    \hbox{\kern-.3em\size{\sf@size}{0pt}\selectfont A}\kern-.15em\TeX}
 \makeatother

  \newcommand{\e}{\eqref}
\newcommand{\ri}{\rightarrow}
\newcommand{\q}{\quad}
\newcommand{\ii}{\infty}

  \newcommand{\ran}{\operatorname{Ran}}
   
   \newcommand{\hess}{\operatorname{Hess}}

\renewcommand\Im{\operatorname{Im}}
\renewcommand\Re{\operatorname{Re}}

\numberwithin{equation}{section}

\newtheorem{lemma}{Lemma}[section]
\newtheorem{theorem}[lemma]{Theorem} 
\newtheorem{corollary}[lemma]{Corollary}
\newtheorem{proposition}[lemma]{Proposition}

\theoremstyle{definition}

\theoremstyle{remark}

\newtheorem{remark}[lemma]{Remark}

\renewcommand{\det}{\operatorname{Det}}
\newcommand{\tr}{\operatorname{Tr}}

  \DeclareMathOperator*{\slim}{s-lim}

\newenvironment{||}{\|}

\def\qqq{\mathrel{\subset\mkern-15mu\lower.38ex\hbox{${\scriptscriptstyle\rightarrow}$}}}
\let\goth\mathfrak
\let\cal\mathcal

\let\Bbb\mathbb

\begin{document}  
 
\title[spectral shift function, trace identities]{Perturbation determinants, the spectral shift function, trace identities, and all that}
\author{ D. R. Yafaev}
\address{ IRMAR, Universit\'{e} de Rennes I\\ Campus de
  Beaulieu, 35042 Rennes Cedex, FRANCE}
\email{yafaev@univ-rennes1.fr}
%\keywords{First order systems, gaps of the spectrum, eigenfunctions,  exponential decay, Dirac and Hill operators}

\begin{abstract}
We discuss applications of the M. G. Kre\u{\i}n theory of the spectral shift function to the multi-dimensional Schr\"odinger operator as well as specific properties of this function, for example, its high-energy asymptotics. Trace identities for  the   Schr\"odinger operator
are derived.
\end{abstract}
\maketitle

%\thispagestyle{empty}

%***********************************************************
\section{Introduction}
%**************************************

 The spectral shift function (SSF)  $\xi(\lambda)=\xi(\lambda; H,H_0)$ is introduced by the relation
\begin{equation}
 \tr  \Bigl( f(H) - f (H_0)\Bigr)=\int_{-\infty}^\infty \xi(\lambda) f^\prime(\lambda)
  d\lambda, 
\label{eq:TF}\end{equation}
called the trace formula. The concept of the
SSF in the perturbation theory  appeared at the beginning of fifties in the
physics literature in the paper \cite{Lif} by I. M. Lifshitz. Its mathematical theory was shortly created by
 M. G. Kre\u{\i}n who proved in \cite{Kr1} relation (\ref{eq:TF})  for a pair of self-adjoint operators $H_{0}$,  $H$
with a trace-class difference $V=H-H_0$ and a wide class of functions $f$. Then it was extended by him in \cite{Kr2}   (see  \cite{Kr3}, for a more complete exposition) to operators $H_{0}$, $H$ with a trace-class difference $R(z)-R_{0}(z)$ of the resolvents. Moreover,  a link  (the Birman-Kre\u{\i}n formula) with the scattering matrix (SM) $S (\lambda)= S (\lambda; H,H_0)$  was found in  \cite{BK} where it was shown that,
for $\lambda$ from a core of the spectrum of $H_{0}$, 
\begin{equation}
 \det  S(\lambda) =\exp\Bigl(-2\pi i \xi(\lambda)\Bigr). 
 \label{eq:BK}\end{equation}

Later   Kre\u{\i}n's theory
was developed in various directions and applied to differential operators. A sufficiently detailed exposition of these developments can be found in survey \cite{BY2} or in book \cite{I} where however applications  to differential operators were almost not discussed. In the present paper which can   to a certain extent be considered as a continuation of  \cite{BY2},  we try to partially fill in this gap.   We study only the Schr\"odinger operator which looks as a sufficiently representative example.

In Section~2 we recall briefly basic results of  Kre\u{\i}n's theory. In Section~3 we discuss its applications to  the Schr\"odinger operator in the space ${\cal H}= L_{2}({\Bbb R}^d)$. Actually,  our goal here is to extend   Kre\u{\i}n's theory to the case of perturbations which are not of trace-class type. Interesting results in this direction can be found in \cite{Kopl1} and \cite{New}. The analogue of the SSF is called the regularized SSF.
Our ``direct" approach turns out to be fruitful even in the case of perturbations   of trace-class type and allows us to reveal   specific properties of the SSF which are not  true for arbitrary self-adjoint operators.  For example, we show that $\xi(\lambda)$ is a continuous function of $\lambda>0$ (in general, it belongs only to the class
 $L_{1}^{(loc)}$), but the proof of this result is curiously complicated for $d > 3$.  
 
 Another goal of this paper is to obtain the high-energy asymptotic expansion of the SSF and of related objects which is the crucial ingredient in derivation of trace identities. Both  these subjects were intensively studied in the literature; see \cite{BF, FaZa} for the one-dimensional case and \cite{Bu, Col} for the case $d  = 3$. We are mainly interested  in the multi-dimensional case where rather stringent assumptions on the potential $v(x)$ (it was required that $v$ belongs to the Schwartz class) were imposed in \cite{Bu, Col}. We note also papers \cite{Rob1, Rob2} where the high-energy asymptotic expansion of the SSF was obtained under natural assumptions on   $v(x)$ in the framework of the microlocal calculus.
 
  In the present paper we suggest a new approach to construction of the asymptotic expansion (AE) of $\xi(\lambda)$ as $\lambda\ri\ii$. Actually, two different problems arise here. The first of them is the proof of the existence of such expansion. The second problem    consists in obtaining   reasonably simple expressions for asymptotic coefficients. We solve the first problem in Section~6  using  the high-energy AE  of the SM $S(\lambda)$ (see Section~5) and the Birman-Kre\u{\i}n formula. This entails also  the AE of the function $\tr ({R^m(z)-R^m_{0}(z)})$, where $m>d/2-1$,  as $|z|\ri\ii$ in the whole complex plane cut along ${\Bbb R}_{+}$. In Section~4 we use another approach to construction of this AE which works only away from the spectrum but gives   reasonably simple expressions for  asymptotic coefficients. This yields also an efficient expression for coefficients in the AE of the SSF.
  Another expression for these coefficients is given in terms of the so called heat invariants. Finally,  trace identities, both of integer and half-integer orders, are derived
  in Section~7.

%***********************************************************
\section{Kre\u{\i}n's theory}
%****

 {\bf 1.}
 Below we use the following standard definitions (see, e.g., \cite{GK}):
 ${\goth S}_{p}$, $p\geq 1$, is the symmetrically normed ideal (with the norm 
 $\|\cdot\|_{p}$) of the algebra of all bounded operators in  a Hilbert space ${\cal H}$; in particular,  ${\goth S}_{1}$ is the ideal of trace operators and ${\goth S}_{2}$ is the ideal of Hilbert-Schmidt operators; the determinant $\det(I+A)$ is well defined for $A \in {\goth S}_{1}$ and possesses standard properties which are basically the same as in the  finite-dimensional case; more generally,  the regularized determinant $\det_{p}(I+A)$, $p=1,2, \ldots$, is well defined for $A \in {\goth S}_{p}$.
 
 Let $H_0$,  $H$ be self-adjoint operators   in  a Hilbert space ${\cal H}$,
 and let $R_0(z)=(H_{0}-z)^{-1}$, $R (z)=(H-z)^{-1}$ be their resolvents.
  The perturbation determinant (PD) for a   pair   $H_0$, $H$   with the difference $V=H-H_{0}\in {\goth S}_{1}$    is defined by the relation
\begin{equation}
D(z)=D_{H/H_0}(z)=\mathrm{Det}\,(I+VR_0(z)) 
\label{eq:PD}\end{equation}
on the set $\rho(H_0)$ of regular points of the  operator $H_{0}$.
 The function $D(z)$ is holomorphic on   $\rho(H_0)$, 
  \begin{equation}
   D (\bar{z})=\overline{D(z)},
   \label{eq:Dp}\end{equation}
    and it has  a zero $z  $ of order $k$ if and only if $z$ is an eigenvalue of multiplicity $k$ of the operator $H$. Moreover,  $D(z)$ satisfies the identity
\begin{equation}
 D^{-1}(z)D^\prime(z)=\mathrm{ Tr}\,(R_0(z)-R(z)),\quad z\in \rho(H_0)\cap\rho(H).
\label{eq:PDR}\end{equation} 
  Since $D(z)\ri 1$ as $|\Im z|\ri\infty$, we can fix    the branch of the function
  $\ln D(z)$ in the half-planes $\pm \Im z > 0$  by the condition $\arg D(z)\ri 0$ as $|\Im z|\ri\infty$. These properties of PD extend to the case $V R_0(z) \in {\goth S}_{1}$.

 The SSF is constructed in terms of the PD in  the following theorem of M. G.  Kre\u{\i}n.

\begin{theorem}\label{Kre}
  For $V\in{\goth S}_1$,  there is the representation 
\begin{equation}
\ln D(z)=\int_{-\infty}^\infty \xi(\lambda)(\lambda-z)^{-1}d\lambda,\quad \Im z\neq 0,
\label{eq:PD1}\end{equation}
where
\begin{equation}
\xi(\lambda)=\pi^{-1}\lim_{\varepsilon\rightarrow +0}\arg
D(\lambda+i\varepsilon).
\label{eq:PD2}\end{equation}
For almost every $($a.e.$)$ $\lambda\in{\Bbb R}$ the limit in $(\ref{eq:PD2})$ exists and
\[
\int_{-\infty}^\infty |\xi(\lambda)|d\lambda \leq \|V\|_1, \q \int_{-\infty}^\infty \xi(\lambda) d\lambda = \mathrm{Tr}\,V.
\]
Moreover, $\pm\xi(\lambda)\geq 0$ if $\pm V\geq 0$.
\end{theorem}

  In a gap of the continuous spectrum  $\xi(\lambda)$ depends on the shift of   eigenvalues of the
operator $H$ relative to   eigenvalues of $H_0$.  

\begin{proposition}\label{Kre1d}
  On component intervals of the set   of common regular points of the operators
$H_0$ and $H$ the SSF $\xi(\lambda)$ assumes constant integral values. If $\lambda$ is an
isolated eigenvalue of finite multiplicity $k_0$ of the operator $H_0$ and $k$ of the operator
$H$, then 
$\xi(\lambda +0)-\xi(\lambda -0)=k_0-k$.
\end{proposition}

The trace formula (\ref{eq:TF}) was justified by M. G. Kre\u{\i}n for a sufficiently wide class of functions. 

\begin{theorem}\label{Kre1}
  Suppose $V\in{\goth S}_1$ and the function $f$ is continuously differentiable while its
derivative admits the representation
\begin{equation}
f^\prime(\lambda)=\int_{-\infty}^\infty \exp(-it\lambda)dm(t),\quad |m|({\Bbb R})<\infty,
\label{eq:Kre}\end{equation}
with a finite $($complex$)$ measure $m$. Then
\begin{equation}
 f(H)-f(H_0)\in{\goth S}_1
\label{eq:Kre1}\end{equation}
and formula $(\ref{eq:TF})$ holds.
\end{theorem}

\begin{remark}\label{BS}
Condition \e{eq:Kre} can be somewhat modified and relaxed (see \cite{BS}, \cite{Pe}).
For example, in the case of bounded operators $H_{0}$, $H$ it suffices to require that
$f^\prime$ is a H\"older continuous function. 
\end{remark}

 \medskip

{\bf 2.}
Now we recall briefly basic notions of 
 scattering theory (see, e.g., \cite{I}, for details).
 Let $E(\cdot)$ be the spectral projection of the operator  $H$, let ${\cal H}^{(a)}$ be  its absolutely continuous subspace and  let
 $  P $ be the orthogonal projection on   ${\cal H}^{(a)}$; the same objects  for the operator $H_{0}$ are endowed with the index $``0"$.
  The wave operator    for a pair of selfadjoint operators $H_0$  and $H$  is defined by  the  relation
\begin{equation}
 W_\pm =W_\pm (H,H_0 )=\slim_{t\rightarrow\pm\infty}\exp(iHt)
\exp(-iH_0t)P_0  
\label{eq:cl4}\end{equation}
 provided   the corresponding strong limit exists.    The wave operator is isometric on ${\cal H}_0^{(a )}$ and  enjoys the intertwining
property
$ W_\pm  E_0(\Lambda)= E (\Lambda) W_\pm  $
where $\Lambda\subset{\Bbb R}$ is an arbitrary Borel set.
In particular, its range $\mathrm{Ran}\: W_\pm \subset {\cal H}^{(a)}$.
The operator $W_\pm $ is said to be complete if    
$  \ran W_\pm={\cal H}^{(a)}$.
 Thus, if the wave operators
 $W_\pm  $ exist and  are complete, then  the absolutely
continuous parts of the   operators $H_0 $ and $H $ are
unitarily equivalent. In this case the scattering operator
\[
{\bf S}= {\bf S} (H,H_0 )= W_+^* (H,H_0 ) W_- (H,H_0 ) 
 \]
 commutes with $H_{0}$ and is unitary on the space ${\cal H}_{0}^{(a)}$.
 
 To define the  SM $S (\lambda; H,H_0 )$, we suppose for simplicity that the spectrum $\sigma_{0}=\sigma (H_{0})$ of the operator $H_{0}$ is absolutely continuous, consists of a finite (or locally finite) union   of closed intervals and has constant multiplicity.
 Let ${\goth h}$ be an auxiliary space whose dimension is equal to this multiplicity.
 Then there exists a unitary mapping ${\cal F}: {\cal H}\ri L_{2}(\sigma_{0}; {\goth h})$ such that the operator ${\cal F}H_{0}{\cal F}^*$ acts as multiplication by the independent variable ($\lambda$) in the space $L_{2}(\sigma_{0}; {\goth h})$. It follows from the relation ${\bf S} H_{0}= H_{0} {\bf S}$ that the operator ${\cal F} {\bf S}  {\cal F}^*$ acts in the space $L_{2}(\sigma_{0}; {\goth h})$ as multiplication by the unitary operator-valued function 
 $S (\lambda)=S (\lambda; H,H_0 ): {\goth h}\ri {\goth h}$ defined for  a.e. $\lambda \in \sigma_{0} $. Of course, this definition fixes the  SM  $S (\lambda  )$ up to a unitary equivalence in the space ${\goth h}$.
 
 The   Kato-Rosenblum theorem asserts that the wave operators    
 $W_\pm (H,H_0 )$ exist and are complete if $V\in {\goth S}_{1}$. 
  A link between the SSF and the SM was found in  \cite{BK}.
  
  \begin{theorem}\label{BK}
  If $V\in{\goth S}_1$, then $S(\lambda)-I  \in {\goth S}_{1}$ and formula \e{eq:BK} holds
  for a.e. $ \lambda\in  \sigma_0$.
\end{theorem}

  Relation  \e{eq:BK} is often used for the definition (up to an integer number) of the SSF on the
absolutely continuous part of the spectrum.  Differentiating it formally, we obtain that
 \begin{equation}
\tr  \Bigl( S^*(\lambda) S ^{\prime}(\lambda)  \Bigr) = -2\pi i \xi^{\prime}(\lambda). 
\label{eq:BKB}\end{equation}
 This relation is sometimes more convenient than  \e{eq:BK} in applications since there exists a simple expression for the trace  of an integral operator with smooth kernel.

  \medskip

{\bf 3.}
 In view of our applications to the Schr\"odinger operator, 
let us consider  generalizations of Theorems~\ref{Kre}, \ref{Kre1} and \ref{BK} specific for the semibounded case.
  With a shift by a constant, it may be achieved that
the operators $H_0+cI$ and $H+cI$ are positive definite. We assume that
\begin{equation}
 R^m(z)-R_0^m(z)\in{\goth S}_1 
\label{eq:4.2.4}\end{equation}
for $z=-c$ and some $m > 0$. 
   It follows from Theorem~\ref{Kre} that the SSF for the pair $h_{0}=(H_0+cI)^{-m}$,
$h=(H+cI)^{-m}$ exists and belongs to the space $L_1({\Bbb R})$. Moreover, according to Proposition~\ref{Kre1d} 
it equals to zero  on ${\Bbb R}_-$ and for sufficiently large $\lambda>0$. For the initial pair
$H_0$, $H$,    the SSF is defined by the equalities
    \begin{equation}  
 \xi(\lambda;H,H_0)=-  \xi((\lambda+c)^{-m};(H+cI)^{-m},(H_0+cI)^{-m}) 
 \label{eq:SSFis}\end{equation}
 for   $\lambda>-c$ and $\xi(\lambda;H,H_0)=0$ for $ \lambda\leq -c$. The last requirement fixes the SSF uniquely.  In terms of the PD $D_{h/h_{0}}(\zeta)$  for the pair $h_{0}, h$, the SSF can be constructed by the formula
 \begin{equation} 
\xi(\lambda; H,H_0 )=\pi^{-1} \arg D_{h/h_{0}} ((\lambda+c+i0)^{-m}). 
\label{eq:5.2.13sm}\end{equation} 
 The class of admissible functions can be obtained from Theorem~\ref{Kre1} (see also Remark~\ref{BS}) by  the change of variables $\mu=(\lambda+c)^{-m}$.
  It follows from the Birman invariance principle that under assumption \e{eq:4.2.4} the wave operators $W_{\pm}(H,H_{0})$ exist and are complete.   In this case the SM are related by the formula
  $S(\lambda;H,H_0)=S^\ast(\mu; h,h _0 )$.  
Thus,   the following result is a direct consequence  of Theorems~\ref{Kre},
  \ref{Kre1} and \ref{BK}. 

\begin{theorem}\label{SB}
Let condition $(\ref{eq:4.2.4})$ hold. Define the SSF $\xi(\lambda )=\xi(\lambda;H,H_{0})$ by    equalities \eqref{eq:SSFis}   for 
  $\lambda>-c$ and $\xi(\lambda)=0$ for $ \lambda\leq -c$.
  It satisfies the condition 
   \begin{equation}
\int_{-\infty}^\infty |\xi(\lambda )|(1+|\lambda|)^{-m-1}d\lambda<\infty  
\label{eq:xi2}\end{equation} 
and is related to the SM by equality $(\ref{eq:BK})$. 
   Suppose that a function $f$ has two locally bounded
derivatives and
  \begin{equation}
 (\lambda^{m+1}f^\prime(\lambda))^\prime= O(\lambda^{-1-\varepsilon}),\quad
\varepsilon>0,\quad
\lambda\rightarrow\infty. 
\label{eq:5.1.4+}\end{equation}
 Then inclusion \e{eq:Kre1} and the trace formula   $(\ref{eq:TF})$ hold.
 \end{theorem}
 
  \begin{corollary}\label{KRYcor}
The trace formula \e{eq:TF} holds  for  all functions $f(\lambda)=(\lambda-z)^{-k}$ where $k\geq m$ and $z\in \rho (H_{0}) \cap \rho (H)$. In particular,
 \begin{equation} 
\tr \Bigl(
R^m (z)-R^m _0(z) \Bigr)=-m \int_{-\infty}^\infty \xi  (\lambda)(\lambda-z)^{-m-1}d\lambda.
\label{eq:ssfd0}\end{equation} 
 \end{corollary}
   
   \begin{remark}\label{GPDe}
If   condition  (\ref{eq:4.2.4})  holds for $m =1$, then
    the generalized PD
  \begin{equation}
  \tilde{D}_{-c}(z):= \det ( I+ (z+c) R(-c) V R_{0}(z))=D_{h/h_0}((z+c)^{-1})
  \label{eq:PDgen}\end{equation}
is correctly defined and satisfies equation \e{eq:PDR}. 
Let us fix the continuous branch of $\arg\tilde{D}_{-c}(z)$ by the condition $\arg\tilde{D}_{-c}(-c)=0$.
 Then it follows from formula \e{eq:5.2.13sm} that
 \begin{equation} 
\xi(\lambda; H,H_0 )=\pi^{-1} \arg \tilde{D}_{-c}(\lambda+i0). 
\label{eq:5.2.13}\end{equation}
\end{remark}

    The SSF has the definite sign for sign-definite perturbations of
the operator  $H_0$. We understand the sign of the perturbation in the sense of quadratic forms.
The following result is  obtained in \cite{Kopl2}  (see also Theorem~8.10.3 of \cite{I}).

\begin{theorem}\label{SBS}
Under assumption $(\ref{eq:4.2.4})$,  $\xi(\lambda )\geq 0$
if $H\geq H_0$ and  $\xi(\lambda )\leq 0$ if $H\leq H_0$.
 \end{theorem}
 
 \medskip

{\bf 4.}
In applications to differential operators it is convenient to
use the concept of  regularized determinant  and to introduce   the
regularized PD
\begin{equation}
 D_p(z)= \mathrm{Det}_p\,(I+ VR_0(z)),\quad
p=2,3,\ldots.
\label{eq:RPD2}\end{equation}
This definition is good for  $VR_0(z)\in{\goth S}_p$. In this case the function $D_{p}(z)$ is holomorphic on the set $\rho(H_0)$.   
    Many properties of   ordinary PD carry over
to regularized PD. Thus, identity  \e{eq:Dp} remains true, and the generalization of (\ref{eq:PDR}) has the form
\begin{equation}
 D_p^{-1}(z)D_p^\prime(z)=-\mathrm{ Tr}\,\Bigl(R(z)-\sum_{k=0}^{p-1}(-1)^k R_0(z)(VR_0(z))^k\Bigr).
\label{eq:PDRp}\end{equation} 

 \medskip

 %%%%%%%%%%%%%%%%%%%%%%%%%%%%%%%%%%%
%%%%%%%%%%%%%%%%%%%%%%%%%%%%%%%%%%%
\section {The  regularized PD and SSF for the Schr\"odinger operator}
 %%%%%%%%%%%%%%%%%%%%%%%%%%%%%%%%%%%%%
%%%%%%%%%%%%%%%%%%%%%%%%%%%%%%%%%%%

{\bf 1.}
Below we consider the pair of self-adjoint operators $H_{0}=-\Delta$,
\begin{equation}
H =-\Delta + v(x), \q v(x)=\overline{v(x)},
\label{eq:S1}\end{equation}
in the space ${\cal H}= L_{2}({\Bbb R}^d)$. It is assumed that the potential $v(x)$ decays sufficiently rapidly at infinity, that is
\begin{equation}
|v(x)| \leq C (1+|x|)^{-\rho},
\label{eq:S2}\end{equation}
where at least $\rho>1$. Then the positive spectrum of the operator $H$ is absolutely continuous, and its negative spectrum consists of eigenvalues
$\lambda_{1}< \lambda_{2}\leq \lambda_{3}\leq\ldots$ counted with their multiplicity (see, e.g.,  \cite{RS}).
 Wave operators \e{eq:cl4} exist and are complete  so that
$\ran W_{\pm} = E ({\Bbb R}_{+}) {\cal H}$ (see, e.g., \cite{Ku3}).

The canonical diagonalization of the operator $H_0$ is realized by the operator
${\cal F}: {\cal H} \ri L_{2 } ({\Bbb R}_{+}; L_{2}({\Bbb S}^{d-1}))$ defined by the equation 
 $  ({\cal F} f)(\lambda;\omega)=(\Gamma_{0}(\lambda)f)(\omega)$ where
 \[
 (\Gamma_{0  } (\lambda) f)(\omega)=2^{-1/2}\lambda^{(d-2)/4}\hat{f}(\lambda^{1/2}\omega), \quad \lambda>0,\quad
\omega\in{\Bbb S}^{d-1}, 
\]
and $\hat{f} $ is the Fourier transform of $f$.   The SM $S(\lambda)=S (\lambda;H,H_0) $ is correctly defined and is a unitary operator in the space
$  L_{2}({\Bbb S}^{d-1})$ for a.e. $\lambda\in {\Bbb R}_{+}$. Moreover, it
   admits the representation
 \begin{equation}
S(\lambda)= I - 2 \pi i ( {\cal  A}_{0}(\lambda) +\tilde{\cal A}(\lambda))
\label{eq:SMst}\end{equation}
 where  
 \begin{equation}
 {\cal A}_0(\lambda)= \Gamma_0(\lambda) V \Gamma_0^\ast(\lambda),
 \q \tilde{\cal A}(\lambda)= -\Gamma_0(\lambda) VR(\lambda+i0)V \Gamma_0^\ast(\lambda).
\label{eq:2.4.71}\end{equation}
These operators are well defined. Indeed, set $\langle x \rangle^{-r}= (1+|x|^2)^{-r/2}$ (we keep the same notation for the operator of multiplication by this function) and choose some $r> 1/2$. Then according to the Sobolev theorem the operator $\Gamma_0(\lambda)\langle x \rangle^{-r}: {\cal H}\ri  L_{2}({\Bbb S}^{d-1})$ is compact and depends continuously on $\lambda > 0$. The limiting absorption principle asserts that the operator-valued  function $\langle x \rangle^{-r}R_{0}(z) \langle x \rangle^{-r}$ is continuous in the complex plane up to the cut along $[0,\ii)$ with the point $z=0$ possibly excluded. Thus, for $\rho> 1$,   operators \e{eq:2.4.71} can be factored in products of bounded operators depending continuously on $\lambda>0$. It follows that the SM $S(\lambda)$ is a continuous (in the operator topology) function of $\lambda>0$.

If estimate \e{eq:S2} is satisfied with $\rho>d$, then  inclusion \e{eq:4.2.4} holds for all
 $ z\not\in   \sigma (H ) $ and any positive integer   $m >d/2 -1$. This result is a straightforward consequence of the resolvent identity  and of the inclusion 
$\langle x \rangle^{-\rho/2} R_0 (z)\in{\goth S}_2$  if $m=1$ but requires some tricks 
if $m >1$  (see \cite{Ytr} or \cite{RS}). Thus,   the   assertion below follows from Theorems~\ref{SB}
and \ref{SBS}.

\begin{theorem}\label{5.1.4}
Suppose that assumption \e{eq:S2} with $\rho>d$ is satisfied.
  Let  $c> -\lambda_{1}$, and let
  a positive integer $ m$ be such that $2(m+1) >d$. Define the SSF by   
    equality $(\ref{eq:SSFis})$  for 
  $\lambda>-c$ and set $\xi(\lambda;H,H_0)=0$ for $ \lambda\leq -c$. 
  Then condition  $(\ref{eq:xi2})$ is satisfied,
\begin{equation}
\xi (\lambda)=0 \quad \mathrm{for}\quad \lambda<\lambda_1, \q
 \xi (\lambda)=-n
  \quad \mathrm{for}\quad \lambda\in(\lambda_n,\lambda_{n+1}),  
\label{eq:5.3.7D2x}\end{equation}
 and   the SSF  is related to the SM by equality $(\ref{eq:BK})$ for $\lambda>0$.
  If a function $f$ has two locally bounded derivatives and satisfies
condition 
$(\ref{eq:5.1.4+})$, then inclusion 
 \e{eq:Kre1} and the trace formula $(\ref{eq:TF})$ $($in particular,  \e{eq:ssfd0}$)$ hold. Moreover, if $v(x)\geq 0$ $(v(x)\leq 0)$, then
$\xi(\lambda)\geq 0$ $($respectively, $\xi(\lambda)\leq 0)$.
 \end{theorem}

\begin{corollary}\label{5.1.4a} 
  Let $d\leq 3$. Then inclusion \e{eq:4.2.4} and
 formula $(\ref{eq:ssfd0})$ hold  for $m=1$, and the SSF can be recovered via  generalized PD $(\ref{eq:PDgen})$  by   relation
$(\ref{eq:5.2.13})$.
 \end{corollary}

\medskip

 {\bf 2.}
 Now we present a direct approach to construction  of the SSF for
  the     multi-dimensional Schr\"odinger operator \e{eq:S1}. This allows us to study some specific properties of $\xi(\lambda)$, such as its continuity   for $\lambda>0$ and   behavior as $\lambda\rightarrow+\infty$ and $\lambda\rightarrow 0$. 
We require here condition (\ref{eq:S2}) for at least $\rho>1$  but not necessarily for $\rho>d$ as in the trace-class approach. Actually, $\rho$ depends on the dimension $d$ of the problem and is different in different assertions. For main results, we assume that 
 $\rho>2$.  In this case the operator $H$ might have only a finite number $ \mathsf{ N}$ of negative eigenvalues.  Thus, the direct approach allows one to introduce the SSF outside of the trace-class framework.
 
In the previous section the SSF was defined by relation (\ref{eq:5.2.13sm}) which is not always convenient.
In the one-dimensional case $VR_{0}(z)\in{\goth S}_{1}$ so that the SSF can be expressed by relation \e{eq:PD2} via the   PD  $D(z)=D_1(z)$ for the original pair $H_0 $, $H$.
To achieve the same for $d>1$, we have to use the regularized PD $D_p(z)$
 defined by formula (\ref{eq:RPD2}).
 The structure of the PD $D_p(z)$ is more complicated for bigger $p$ and usually some property of $D_p(z)$ entails the same property of  PD for bigger $p$. Therefore it is natural to choose the smallest   $p$ such that $VR_{0}(z)\in{\goth S}_{p}$.

The following elementary assertion   shows that   properties of regularized PD  are essentially the same as in the trace-class case.

\begin{proposition}\label{PDG} 
Suppose that   $\rho>1$ if $d=1, 2$,  $\rho>3/2$ if $d=3$ 
and $\rho\geq 2$ if $d\geq 4$. Let $p = 1$ for $d=1$, $p =2$ for $d=2, 3$ and  $p> d/2$ for $d \geq 4$.
 Then the function $D_p(z)$ is analytic
in the complex plane cut along the positive half-axis $[0,\infty)$ and satisfies identities
 \e{eq:Dp} and \e{eq:PDRp}.  For  $z\not\in[0,\infty)$, this function may have zeros  on  the negative half-axis only,
where they coincide with   eigenvalues of the operator $H$; the order of a zero of $D_p(z)$ equals
the multiplicity of the corresponding  eigenvalue. 
 \end{proposition}

Let $r=\rho/2$,  $G= \langle x \rangle^{-r}$ and let ${\cal V}$ be the multiplication operator by the bounded function
$\langle x \rangle^{\rho} v(x)$.
Since non-zero eigenvalues of the operators $V R_0(z)$ and  $G R_{0}(z) G{\cal V}$  are the same, we have that
\[
D_p(z)=  \det_p (I+  G R_{0}(z) G   {\cal V} ).
\]
This representation is more convenient than (\ref{eq:RPD2}) because the operator-valued  function
 $ G R_{0}(z)  G $  is continuous   up to the cut even in the classes ${\goth S}_{p}$ for sufficiently large $p$. In the following assertion we take also into account that if $D_p(\lambda\pm i0)=0$ for $\lambda>0$, then $-1$ is an eigenvalue of  the compact
operator $  G R_0(\lambda\pm i0)G{\cal V}$  and hence  
$\lambda$ is an eigenvalue of $H$. However according to the Kato theorem (see, e.g., \cite{RS}) this operator does not have positive eigenvalues.
   
\begin{proposition}\label{PDG1} 
Let    $\rho> 1$ if $d=1$, $\rho>3/2$ if $d=2$  and  $\rho>2$ if $d\geq 3$. Put $p=1$ for $d=1$, $p=2$ for $d=2, 3$  and  $p= d$ for $d \geq 4$.
 Then the function $D_p(z)$   is continuous
 up to the cut along $[0,\ii)$, with the point
$z=0$ possibly excluded. The  function $D_p(\lambda\pm i0)$ does not have zeros for $\lambda> 0$.
   \end{proposition} 
    
 As far as the high-energy behavior of the PD is concerned, we consider only the case $d=3$.  Let us use the
identity  
\begin{equation} 
 \det_2\,(I+G R_0(z)G {\cal V})=\det_3\,(I+ G R_0(z)G {\cal V})
\exp (-2^{-1} \tr (VR_0(z))^2).
\label{eq:5.3.3}\end{equation}
Remark that
\[ 
\| G R_0(z) G \|_3^3\leq 
 \| G R_0(z) G\|
\,  \| G R_0(z)G\|_2^2. 
\]
The first factor in the right-hand side tends to zero as $|z|\rightarrow\infty$, and it follows from the explicit formula for the integral kernel of the operator $R_{0}(z)$ that the second
factor is uniformly bounded. Therefore  the first factor in the right-hand side of (\ref{eq:5.3.3}) tends to $1$ as $|z|\rightarrow\infty$. The same is true for the second factor because
\[
\tr (VR_0(z))^2 =(4\pi)^{-2}\int_{{\Bbb R}^3}\int_{{\Bbb
R}^3} v(x) v(y) |x-y|^{-2}e^{2i \sqrt{z} |x-y|} dxdy,
\]
and the integral tends to zero according to the Riemann-Lebesgue lemma. 
More generally, we  have 

\begin{proposition}\label{5.3.1} 
Under the assumptions of Proposition~$\ref{PDG1}$  
\begin{equation} 
\lim_{|z|\rightarrow\infty} D_p(z)=1
\label{eq:5.3.2}\end{equation}
uniformly in $\arg z\in [0,2\pi]$.
\end{proposition}

The PD $D_p(z)$ is singular at the point $z = 0$ for $d=1$ and $d=2$ only.

\begin{proposition}\label{PDG2} 
Let  $d \geq 3$,  $\rho> 2$, $p=2$ for $d=3$   and  $p= d$ for $d \geq 4$.
 Then the PD $D_p(0)$ is correctly defined and $ \| D_p(z)-   D_p(0)\| \ri 0$ as 
 $|z|\rightarrow 0$.
\end{proposition}
 
Clearly, $D_p(0)=0$ if and only if $-1$ is an eigenvalue of the compact operator 
$   G R_0(0)G {\cal V}$. In this case one says that the operator $H$ has a zero-energy resonance (in particular, $H$ might have zero eigenvalue).

We avoid a study of the PD $D_p (z)$ as $ | z | \rightarrow 0$ imposing a mild a priori assumption
    \begin{equation} 
 \lim_{|z|\rightarrow 0}|z|^\alpha  \ln D_{p} (z) =0
\label{eq:7.3.15}\end{equation}
for a suitable $\alpha>0$. Of course
  in the    case $d\geq 3$   condition (\ref{eq:7.3.15}) is satisfied for all $\alpha>0$  if $D_{p}(0)\neq 0$.   More generally,    condition    (\ref{eq:7.3.15}) 
 can be deduced from the low-energy expansion (as $z\ri 0$) of $R(z)$  which requires however  a somewhat  stronger assumption  than   $\rho>2$.

\medskip

{\bf 3.}
We suppose that   $\arg D_p(z)$ is a continuous  function of $z$ in the complex plane cut along $[\lambda_{1},\ii)$  which is possible because $  D_p(z)\neq 0$ there.
According to (\ref{eq:5.3.2}) we then  fix its branch  by the condition 
$   \arg  D_p(z)\ri  0$ as $  |z|\rightarrow \infty$.

 Now we are in a position to construct the regularized SSF $\xi_p$.

\begin{theorem}\label{5.3.2} 
 Let condition $(\ref{eq:S2}) $ where $\rho>2$  hold,  and let  $p= 1$ for $d=1$,
  $p = 2$ for $d=2, 3$  and  $p =d$ for $d \geq 4$. 
  Assume that 
\e{eq:7.3.15} is true for $\alpha=1$.
        Define the regularized SSF $\xi_p$ by   equality \e{eq:PD2} in terms of the regularized PD $D_{p}(z)$.
Then        representation \e{eq:PD1}  for $ \ln D_p(z)$  holds with the function $\xi_p(\lambda)$.
The  function $\xi_p$ is
determined by   equalities   \e{eq:5.3.7D2x}  for  $\lambda<0$.  
It is  continuous for $\lambda>0$,
$  \xi_{p} (\lambda)=o(1)$ as $ \lambda\ri \infty$,  
 $\xi_{p} (\lambda)=o(\lambda^{-1})$ as $\lambda\ri 0$
 and    the integrals of $\xi_{p} (\lambda) \lambda^{-1}$ and $\xi_{p} (\lambda)$
are    convergent $($but  not necessarily absolutely$)$  at the points $\lambda=\infty$ and $\lambda=0$, respectively.
\end{theorem}

Indeed, the function    $ \ln   D_{p} (z)$ 
is   analytic in the complex plane cut along $[\lambda_1 ,\infty)$
   and is continuous up to the cut with
exception of the points $\lambda_1 ,\ldots, \lambda_\mathsf{N} $ and, possibly, zero.
Let us consider in the complex plane the closed contour $\Gamma_{R,\varepsilon}$ which consists of  the intervals $(\lambda_1,R+i0)$ and $(R-i0, \lambda_1)$ lying on the upper and lower edges of the cut and of the circle
$C_R$ of radius $R$   passed   in the counterclockwise direction. 
   Moreover, we bypass every
point $\lambda_j$, $j=1 ,\ldots, \mathsf{N}$, and the point $ 0$  by      semicircles $C^\pm_\varepsilon(\lambda_j)$ and $C^\pm_\varepsilon(0)$ of radius
$\varepsilon$. By virtue of the Cauchy theorem, for an arbitrary complex $z$,
 a  sufficiently small $\varepsilon$ and a sufficiently large $R$
\begin{equation}  
\ln D_{p}(z)=(2\pi i)^{-1}\int_{\Gamma_{R,\varepsilon}} 
\ln D_{p} (\zeta)(\zeta -z)^{-1} d \zeta. 
\label{eq:5.2.23}\end{equation}
According to \e{eq:5.3.2}, the integral over $C_R$ tends to zero as
$R\rightarrow\infty$. Since the function $D (z)$ has only   zeros of finite order at the points
$\lambda_1,\ldots, \lambda_\mathsf{N}$, the integrals over $C^\pm_\varepsilon(\lambda_j)$ tend  to zero as
$\varepsilon\rightarrow 0$. The integrals over $C^\pm_\varepsilon(0)$ tend to zero as
$\varepsilon\rightarrow 0$ according to condition (\ref{eq:7.3.15}). Hence representation \e{eq:PD1} for $D_{p}(z)$ with   function $\xi_{p}(\lambda)$ defined by \e{eq:PD2} follows from equations
  (\ref{eq:Dp}) and  (\ref{eq:5.2.23}).
 
Let us now obtain the trace formula.

\begin{theorem}\label{5.3.2tr}   
 Let $f$ be a  bounded rational function with non-real poles.
Then under the assumptions of Theorem~$\ref{5.3.2}$
 \begin{equation} 
\tr \Bigl(
f (H)-\sum_{k=0}^{p-1} \frac{1}{k!}  \frac{d^k  f (H_{0}+\varepsilon V)}{d\varepsilon^k}\Big|_{\varepsilon=0} \Bigr)= \int_{-\infty}^\infty\xi_p (\lambda)
f^{\prime}(\lambda)d\lambda
\label{eq:5.3.ph}\end{equation}
and, in particular, 
 \begin{equation} 
\tr \Bigl(
R(z)-\sum_{k=0}^{p-1}(-1)^k R_0(z)(VR_0(z))^k \Bigr)=-\int_{-\infty}^\infty\xi_p(\lambda)(\lambda-z)^{-2}d\lambda.
\label{eq:5.3.10}\end{equation}
\end{theorem}

Indeed, differentiating representation (\ref{eq:PD1}) for $D_{p}(z)$ and taking into account   formula (\ref{eq:PDRp}),
we obtain first \e{eq:5.3.10} which yields \e{eq:5.3.ph} for $ f (\lambda)=(\lambda-z)^{-1}$, $\Im z\neq 0$.  Then differentiating \e{eq:5.3.10}, we extend  formula  \e{eq:5.3.ph} to functions
  $f (\lambda)=(\lambda-z)^{-m}$ where $m=2,3,\ldots$. Clearly,   formula \e{eq:5.3.ph} remains  true for linear combinations of such functions.

The trace formula  (\ref{eq:TF})  makes of course no sense if $\rho\leq d$, and  
representation \e{eq:5.3.ph}  can be regarded as its regularization.  If $d=1$, then $\xi=\xi_{1}$ and Theorem~\ref{5.3.2tr} is contained in Theorem~\ref{SB}.

Applying the argument principle to the function  $D_{p}(z) $, we obtain the trace identity of zero order.

\begin{theorem}\label{5.3Lev} 
Under the assumptions of Proposition~\ref{PDG2} suppose that $D_{p}(0)\neq 0$. Then
\begin{equation}  
\arg D_p(\infty)-\arg D_p(0)=\pi \mathsf{N}, 
\label{eq:5.3.11}\end{equation}
where $\mathsf{N}$ is the total number of negative eigenvalues of the
operator $H$.
\end{theorem}

We emphasize that according to definition \e{eq:PD2}
under the assumption $D_{p}(0)\neq 0$ the limit
$\xi_{p}(+0) $ exists. Formula \e{eq:5.3.11} (known as the Levinson formula) means that the regularized SSF 
$\xi_{p}(\lambda) $ is continuous at the point $\lambda=0$. In the cases $d=1$ and $d=2$ the singularity of the PD $D_{p}(z) $ at the point $z=0$ should be taken into account.

For positive $\lambda$, the regularized SSF is related to the SM $S(\lambda)$.
For simplicity we consider the case $p=2$ only.  Let
  the operator ${\cal A}_0(\lambda)$ (known as the first Born approximation to the SM)  be defined by formula  \e{eq:2.4.71}. It is self-adjoint, does
not depend on the resolvent of
$H$ and is constructed directly in terms of the Fourier transform of the  potential $v$.  We set
 \[ 
 \nu(\lambda)=  \tr\Bigl(S(\lambda)-I +2 \pi i {\cal A}_{0}(\lambda) \Bigr)
=- 2 \pi i \tr \tilde{\cal
A}(\lambda).
\]
 The next result can be considered as a modification of the  Birman-Kre\u{\i}n formula $(\ref{eq:BK}) $.  Its proof     relies on the identity    
 \begin{eqnarray*} 
 \det_2\,\Bigl(I-(V-V R(z)V) (R_0(z) -R_0(\bar{z}))\Bigr)
\nonumber\\
=D_2(\bar{z}) D_2(z)^{-1}\exp \Bigl(-\tr (V R (z) V  (R_0(z) -R_0(\bar{z})))\Bigr)
\end{eqnarray*}
where we pass   to the limit $z\rightarrow\lambda+i0$.

\begin{theorem}\label{5.3SMSSF} 
Let condition  \e{eq:S2} be satisfied for  $\rho>3/2$
if  $d=2$ and for $\rho>2$ if $d=3$. Then  the function $\xi_2(\lambda) $
is related to the SM $S(\lambda)$ by the formula
\[
 \det_{2}  S(\lambda)=e^{-2\pi i \xi_2(\lambda)- \nu (\lambda)},\q
 \lambda>0. 
 \]
\end{theorem}

\medskip

 {\bf 4.}
Suppose now that $v$ satisfies condition \e{eq:S2} with $\rho>d$. Then inclusion (\ref{eq:4.2.4})  holds for $2(m+1)> d$, and  the usual SSF $\xi(\lambda)$ is    defined by relation (\ref{eq:SSFis}).
We  will find a relation between  $\xi$ and the regularized SSF $\xi_p$ which
 will allow us to obtain a new  information on the SSF $\xi$ supplementing the results of  subs.~1.

Let first $d\leq  3$ so that $m=1$. The relation between   the generalized PD
 $\tilde{D} (z)=\tilde{D}_{-c}(z)$, $-c < \lambda_{1}$, defined by  equality
(\ref{eq:PDgen}) and  the regularized PD
 $D_{2} (z)$ can easily be obtained by comparison of equations (\ref{eq:PDR}) for $\tilde{D} (z)$ and (\ref{eq:PDRp}) for $p=2$.

 \begin{proposition}\label{delta} 
  If $d=2$ or $d= 3$ and  condition \e{eq:S2} is satisfied with $\rho>d$, then 
  \begin{equation}
\tilde{D}(z)=D_{2}^{-1} (-c)D_{2} (z)\exp\Bigl( (4\pi)^{-1}\int_{{\Bbb R}^d} v(x)dx \, b_{d}(z)\Bigr),
\label{eq:PDRarg}\end{equation}
   where $b_{2}(z)=-\ln (-z/c)$, $b_{3}(z)=c^{1/2} -(-z)^{1/2}$ and $\arg(-z)\in (-\pi,\pi)$.
  \end{proposition} 
   
   Taking the arguments of both sides of \e{eq:PDRarg} and
 passing   to the limit $z\ri\lambda+i0$,
  we obtain the relation between the SSF $\xi  (\lambda)$ which can be defined by formula (\ref{eq:5.2.13}) and $\xi_2 (\lambda)$.
  Then we use  Theorems~\ref{5.3.2}   and \ref{5.3Lev}.

\begin{theorem}\label{5.3.3} 
 Let the assumptions of Proposition~\ref{delta}  be satisfied.
 Then for $\lambda>0$ the SSF $\xi(\lambda) $ admits
the representation
 \[ 
\xi(\lambda)=\xi_2(\lambda )+  2^{-1}(2\pi)^{-d}|{\Bbb S}^{d-1}|
  \int_{{\Bbb R}^d} v (x)dx \;   \lambda^{(d-2)/2}.
  \]
 The  function $\xi(\lambda)$    is continuous for $\lambda> 0$   and   
\begin{equation}
\xi(\lambda)=2^{-1}(2\pi)^{-d}|{\Bbb S}^{d-1}|
  \int_{{\Bbb R}^d} v (x)dx \; \lambda^{(d-2)/2}+ o(1), \q \lambda\rightarrow \infty. 
  \label{eq:ssfhe}\end{equation}
 If $d=3$ and the operator $H$ does not have
the zero-energy resonance, then there exists the limit $\xi(+0)= - \pi \mathsf{N}$   so that
 the SSF is continuous at  the point $\lambda = 0$.
\end{theorem}
  
  If $d=1$, then $D_{2}(z)$ is related to the usual PD $D (z)$ by a formula similar to \e{eq:PDRarg}. This implies that the SSF $\xi(\lambda)$ is continuous for $\lambda> 0$   and satisfies relation \e{eq:ssfhe} where the remainder is $o(\lambda^{-1/2})$.

  The case $d\geq 4$ is considerably more difficult because the SSF $\xi(\lambda)$ is
 expressed only  in terms of the PD for the pair $(H_{0}+cI)^{-m}$,  $(H+cI)^{-m}$ which is not directly related to $D_{p}(z)$ if $m> d/2-1\geq 1$.  To check   that $\xi$ is a continuous function of $\lambda> 0$,   we compare   formulas (\ref{eq:ssfd0})  and
\e{eq:5.3.10} and use that the function
   $ \tr \Bigl( d^{m-2 }      (R_0(z) (VR_0(z))^k  ) / dz^{m-2} \Bigr)$, $k\geq 2$, 
     is continuous up to the cut along $(0,\infty)$. This implies that the integral of $(\xi(\lambda) - \xi_p(\lambda))(\lambda-z)^{-m}$ over $\lambda$  is also a continuous function of $z$ up to this cut. Then 
 standard results on the Cauchy type integrals show that
    the function $\xi(\lambda) - \xi_p(\lambda) $  belongs   to the class $C^{m-1}$  for $\lambda> 0$.    Now the continuity  of   $\xi(\lambda)$ follows from Theorem~\ref{5.3.2}. Let us formulate this result.
    
   \begin{theorem}\label{5.HD} 
 Let   condition  \e{eq:S2} with $\rho>d$   be satisfied.  
Then  the SSF $\xi(\lambda) $      is continuous for $\lambda> 0$.
\end{theorem}

Using representation \e{eq:SMst} it is easy to show that under the assumption $\rho>d$
the SM $S(\lambda)$ depends continuously on $\lambda>0$   in the trace-class topology and hence  $\det S(\lambda)$ is also a continuous function of $\lambda>0$. Nevertheless the continuity of the SSF $\xi(\lambda)$ does not follow from formula (\ref{eq:BK}) because   its integer jumps are not a priori excluded.

The  leading term of the high-energy asymptotics of the SSF is given 
by the first term in the right-hand side of
 \e{eq:ssfhe}  for all $d\geq 4$. This result is proven in Section~6, together with  a complete AE of $\xi(\lambda)$ as $\lambda\ri\infty$, where more stringent assumptions on $v(x)$ are imposed.

 %%%%%%%%%%%%%%%%%%%%%%%%%%%%%%%%%%%
%%%%%%%%%%%%%%%%%%%%%%%%%%%%%%%%%%%
\section {The Green function for large values of the spectral parameter}
 %%%%%%%%%%%%%%%%%%%%%%%%%%%%%%%%%%%%%
%%%%%%%%%%%%%%%%%%%%%%%%%%%%%%%%%%%

In subs.~1 we  construct an AE   of the integral kernel $R (x,x^\prime ; z)$ of the resolvent (of the Green function) as $|z|\ri\infty$ provided $z\in\Pi_{\theta}$ where $\arg z\in (\theta, 2\pi-\theta)$. This method works only for $\theta>0$ but gives explicit expressions for the coefficients of this expansion and does not require any decay of $v(x)$ at infinity. Under some  decay assumptions this method yields also   an AE of   $\tr(R^m(z)-R_{0}^m (z))$.
In subs.~2 we obtain a local   AE of the  parabolic Green function (heat kernel)  $G (x,x^\prime ;t)$ as $t\ri 0$ which requires almost no assumptions on $v(x)$ at infinity. The Laplace transform  relates these expansions with the AE of    $R (x,x^\prime ; z)$ as $|z| \ri \ii$. In subs.~3  these results   are used to enhance the results on the local AE of    $R (x,x^\prime ; z)$. On the contrary,  the   AE of  $\tr(R^m(z)-R_{0}^m (z))$ is used to derive an AE of  
$\tr ( e^{- H  t} -e^{- H_{0}  t}) $  as $t\ri 0$.

\medskip

{\bf 1.} 
We proceed from a modification of the iterated resolvent identity
which is a special case of the non-commutative Taylor formula of \cite{Kant}. We formulate it for the Schr\"odinger
 operator, but   actually this identity has  abstract nature.
 
  \begin{proposition}\label{AsExpTRK}
  Suppose that $v\in C^\infty({\Bbb R}^d)$ and that $v$ as well as all its derivatives are bounded functions.  Define the operators $X_{n}$ by the recurrent relations 
    \begin{equation} 
     X_{0}=I \q \mathrm{and}  \q X_{n+1}=X_{n}H_{0}- H X_{n}
\label{eq:Kant}\end{equation}
so that
 \[
 X_{n}=\sum_{k=0}^n (-1)^k
    \begin{pmatrix}
  n\\ k
    \end{pmatrix}
     H^k H_{0}^{n-k}.
 \]
 Then, for all $N\geq 0$,
   \begin{equation} 
 R(z) =  \sum_{n=0}^N X_{n}R_{0}^{n+1}(z)+ R(z) X_{N+1}R_{0}^{N+1}(z) .
\label{eq:AsExpTRK}\end{equation}
   \end{proposition}
 
  The operators $X_{n}$   can   easily be computed.  For example, $X_{1}=-v$, 
     $
     X_{2}=-2 \langle (\nabla v) , \nabla\rangle - (\Delta v) +v^2
     $ 
     and
      \[
     X_{3}= -4  \langle\hess v \, \nabla,   \nabla\rangle   +6v \langle ( \nabla  v) , \nabla\rangle
   - (\Delta^2 v)  + 2 |\nabla v |^2 +3v (\Delta v) - v^3.
     \]
         
      \begin{proposition}\label{Kant}
    The $X_{n}$ is a differential operator of order $n-1$ so that
     \begin{equation} 
 X_{n} =  \sum_{|\alpha|\leq n-1}p_{\alpha,n}\partial^\alpha,
     \q p_{\alpha,n}(x)=\overline{p_{\alpha,n}(x)}, \q \alpha=\{\alpha_{1},\ldots, \alpha_{d}\}.
\label{eq:AsExpTRK1}\end{equation}
Under the assumptions of  Proposition~\ref{AsExpTRK}
  all coefficients $p_{\alpha,n}$       as well as all its derivatives  are bounded functions.
           Moreover,  if $v$ satisfies estimates
           \begin{equation}
\vert \partial^\kappa v(x)\vert \leq C_\kappa (1+\vert x\vert)^{-\rho -\vert \kappa\vert},\quad \rho >0,  
\label{eq:1.4LR}\end{equation}
 for all multi-indices $  \kappa $ , then $(\partial^\kappaÊp_{\alpha,n})(x)= O( |x|^{-\rho-|\kappa|-(n-1)\epsilon})$, $\epsilon=\min\{\rho,1\}$,
   for all $n$, $\alpha$ and $\kappa$ as $|x|\ri\ii$.
 \end{proposition}
 
Both Propositions~\ref{AsExpTRK}  and \ref{Kant}   can    be verified by induction in $n$ quite straightforwardly.

 Let us use identity \e{eq:AsExpTRK} to obtain an AE of the    integral kernel $R(x,x^\prime;z)$ of $R(z)$ as $|z|\ri\infty$, $z\in \Pi_{\theta}$. According to \e{eq:AsExpTRK1}   the operator $T_{n}(z)= X_{n}R_{0}^{n+1}(z)$ for $n\geq 1$ has integral kernel
    \begin{equation} 
T_{n}(x,x^\prime;z)  = (2\pi)^{-d} \sum_{|\alpha|\leq n-1} i^{|\alpha|} p_{\alpha,n}(x)
      \int_{{\Bbb R}^d} e^{i\langle x- x^\prime,\xi \rangle}   \xi^\alpha(|\xi|^2-z)^{-n-1} d\xi.  
\label{eq:AsyRes}\end{equation}
 It follows  that
      \[ 
| T_{n}(x,x^\prime;z) |\leq C  \sum_{|\alpha|\leq n-1}  | p_{\alpha,n}(x)|
    |z|^{- (n+3-d)/2} 
    \] 
    if $n+3> d$. Of course formula \e{eq:AsyRes} can be rewritten as
       \begin{equation} 
T_{n}(x,x^\prime;z)  = n !^{-1} \sum_{|\alpha|\leq n-1}  p_{\alpha,n}(x)
    \partial_{x}^\alpha R_{0}^{(n)} (x-x^\prime;z)
\label{eq:AsyResR}\end{equation}
    where $R_{0}^{(n)}(x-x^\prime;z)=\partial_{z}^{n} R_{0} (x-x^\prime;z)$ and
     $ R_{0} (x-x^\prime;z)$ is integral  kernel of the resolvent $R_{0}(z)$.
     The estimate of the remainder $ R(z) X_{N+1}R_{0}^{N+1}(z)$ in 
     \e{eq:AsExpTRK} relies only on the trivial bound $\| R(z)\|= O (|z|^{-1})$ and explicit formula \e{eq:AsyRes}.    Thus,  identity  \e{eq:AsExpTRK} yields the following result.
              
               \begin{theorem}\label{AsyRes}
               Suppose that $v\in C^\infty({\Bbb R}^d)$ and that $v$ as well as all its derivatives are bounded functions.   Then, for all sufficiently large $N$  $(N\geq d-3)$, the  asymptotic relation 
       \begin{equation} 
R(x,x^\prime;z) = \sum_{n=0}^N  T_{n}(x,x^\prime;z)+ O (| z|^{-(N-d)/2-2}) 
       \label{eq:AsyRes1}\end{equation}
        is valid as $|z|\ri \infty$, $z\in \Pi_{\theta}$. The estimate of the remainder here  is  uniform  with respect to $x,x^\prime\in  {\Bbb R}^d$. Moreover, relation \e{eq:AsyRes1} can be infinitely differentiated in  $x$, $x^\prime$ $($then $N$ increases$)$ and $z$.
\end{theorem}

 It follows from \e{eq:AsyResR} that the functions $T_{n}(x,x^\prime;z)$ decay exponentially as $|x - x^\prime| |z|^{1/2}\ri \infty$, $z\in \Pi_{\theta}$. Therefore,
although valid for all $x,x^\prime\in  {\Bbb R}^d$, AE  \e{eq:AsyRes1}
is of interest in the region $|x - x^\prime| =O (| z|^{-1/2})$ only.  
Observe that the functions $T_{n}(x,x^\prime;z)$ are singular on the diagonal $x = x^\prime$ but are getting smoother as $n$ increases. In particular, these functions are continuous if $n >  d- 3$ and $n\geq 1$ ($T_0(x,x^\prime;z)$ is continuous for $d=1$ only). Thus, \e{eq:AsyRes1} yields the expansion of the Green function both for large $|z|$ and in smoothness. Note also that diagonal singularities   of the functions $T_n (x, x^\prime; z)$ are getting weaker after differentiations with respect to $z$.

 Let us  discuss  expansion \e{eq:AsyRes1} for the case $x=x^\prime$.
 If $d=1$, then we can directly set $x = x^\prime$ in  \e{eq:AsyRes1}. In the case $d\geq 2$ we previously differentiate  \e{eq:AsyRes} and  \e{eq:AsyRes1} $(m-1)$-times  with respect to $z$.
 If $2(m+1)>d$, then  integral kernel of the operator  $R^{(m-1)} (  z) -R^{(m-1)}_{0} (  z)$ is continuous. Therefore   setting $x=x^\prime$ and  collecting together  terms of the same power of $z$, we obtain
  
   \begin{corollary}\label{Parab4}
  Let $2(m+1)>d$. Set
  $  c_{\alpha,k}=    \int_{{\Bbb R}^d} \xi^{2\alpha} (|\xi|^2+1)^{-k} d\xi$
  where $k> |\alpha| +d/2$.
   Then
  \begin{equation}
\Big (R^{(m-1)}  (x,x^\prime ;z)- R_{0}^{(m-1)}  (x,x^\prime ;z)\Big)\Big|_{x=x^\prime}=    (m-1)! \sum_{n=0}^\infty  
        r_{n}^{(m)} (x) (-z)^{d/2 -m-n -1 }  
\label{eq:loc2a}\end{equation}
where $($the definition below makes sense for $n+m+1> d/2)$
\begin{equation}
 r_{n}^{(m)} (x)= (2\pi)^{-d}  \sum_{|\alpha|\leq n} (-1)^{|\alpha|}
   \begin{pmatrix}
   |\alpha| +n+m\\
   m-1
     \end{pmatrix}
 c_{\alpha,  |\alpha| +n+m+1 }  p_{2\alpha,n +|\alpha| +1}(x). 
\label{eq:AsExpTloc}\end{equation}
 \end{corollary} 
 
     For potentials decaying at infinity,   AE \e{eq:loc2a} can be integrated over $x\in {\Bbb R}^d$ because 
  the derivative of order $m-1$ of the remainder in \e{eq:AsExpTRK} can be estimated in the trace-class norm.     This yields the following result. 
   
   \begin{theorem}\label{AsExpTRMP}
    Let assumption $(\ref{eq:1.4LR})$ where $\rho>d$ hold, and let $2(m+1)>d$.
      Then  the  AE as $|z|\ri\infty$ 
        \begin{equation} 
         \tr(R^m(z)-R_{0}^m (z)) =   \sum_{n=0}^\infty  
       \mathbf{r}_{n}^{(m)}  (-z)^{d/2 -m-n -1 }     
\label{eq:AsExpTRK7}\end{equation} 
is valid for  $z\in\Pi_{\theta}$.   The  real coefficients $ \mathbf{r}_{n}^{(m)}$ are defined by the formula
   \begin{equation}
 \mathbf{r}_{n}^{(m)} =   \int_{{\Bbb R}^d} r_{n}^{(m)} (x) dx.
\label{eq:AsExpTRK6}\end{equation} 
\end{theorem}

\begin{remark}\label{AsExpTRMPr}
   Instead of taking derivatives, we can remove from $R(z)$ several (instead of one as in  \e{eq:AsExpTRK7}) terms of its expansion \e{eq:AsExpTRK}. This allows us   to relax also the condition $\rho>d$. For example, if $d< 6$ and (\ref{eq:1.4LR}) is satisfied  for
    $\rho>d/2$, then 
     \begin{equation}
         \tr(R(z)-R_{0} (z) +R_{0} (z) V R_{0} (z)) =   \sum_{n=1}^\infty  
       \mathbf{r}_{n}^{(1)}  (-z)^{d/2 -n - 2 }.     
\label{eq:AsExpTRK7r}\end{equation}
\end{remark}

\medskip

{\bf 2.}
The  parabolic Green function
 $G (x,x^\prime ;t)$ is   integral kernel of the operator
  $\exp(- H t)$, $t>0$. It satisfies  the parabolic equation
  \begin{equation}
\partial G (x,x^\prime ;t)/ \partial t=\Delta ÊG  (x,x^\prime;t) - v(x) G(x,x^\prime;t),\q \Delta=\Delta_{x}. 
\label{eq:par2}\end{equation}
We seek (cf. \cite{BabRap}) its approximate solution   in the form
   \begin{equation}
G_{N} (x,x^\prime ;t)=G_{0}(x,x^\prime;t) \sum_{n=0}^N g_{n} (x,x^\prime  ) t^n  
\label{eq:par1}\end{equation}
where    $G_{0}(x,x^\prime;t)=G_{0}(x - x^\prime;t)$ is integral  kernel of the operator $\exp(- H_{0} t)$  and $ g_{0} (x,x^\prime  )=1$. Let us plug expression \e{eq:par1} into
 \e{eq:par2}, use the equation $\partial G_{0}   / \partial t= Ê\Delta G_{0}$ and divide by the common factor $G_{0}$.
  Requiring then that the coefficients at $t^n$, $n=0,1,\ldots, N-1$, vanish,  we find recurrent equations for the functions $g_{n+1} (x,x^\prime  )$. For a fixed $x^\prime$, this yields an ordinary differential equation  $(n+1) g_{n+1}   + r \partial g_{n+1} / \partial r = \DeltaÊg_{n }   - v  g_{n }$ in the variable $r=|x-x^\prime|$.
  Solving it under the assumption that
$g_{n +1} (x,x   )$ is finite, we obtain the formula
\begin{equation}
  g_{n+1} (x,x^\prime  )   = \int_{0}^1   \sigma^nÊ(\Delta g_{n }  - v  g_{n }  ) (x^\prime+\sigma (x-x^\prime),x^\prime  )d\sigma. 
\label{eq:par6}\end{equation}
The construction above does not require any assumptions on $v(x)$ at infinity. Our justification of the local AE of   $G (x,x^\prime ;t)$ as $t\ri 0$ requires only very mild assumptions. Using representation \e{eq:par6}, one easily proves
 the following assertion by induction.

  \begin{lemma}\label{Parab}
      Suppose that $v\in C^\infty({\Bbb R}^d)$ and $(\partial^\kappa v)(x)= O(e^{\gamma |x|^2})$ as $|x|\ri\infty$ for some $\gamma\geq 0$ and all $\kappa$. Then  
       $ |   g_{n}  (x,x^\prime   ) |
    \leq C_{n }    e^{ (n\gamma+\varepsilon)  |x - x^\prime |^2} $ for all  $ \varepsilon>0$     
       uniformly   with respect to $x^\prime$ from a compact subset of ${\Bbb R}^d$. Moreover these estimates can be infinitely differentiated in the variables
$x$ and $x^\prime$.
 \end{lemma}

Let us set $Q_{N}= G_{0} ( - \Delta Êg_{N } + v  g_{N }  )t^N$. It follows from \e{eq:par2} and a  similar equation for the function $G_{N}$ that  the difference
  $ F_{N}  = G  -G_{N} $ admits    the representation
\begin{equation}
F_{N} (x,x^\prime ;t) = -\int_{0}^t (e^{-(t-s)H} Q_{N}(\cdot,x^\prime ;s))(x)ds 
\label{eq:par7f}\end{equation} 
where a point $x^\prime$  is fixed. Combining estimate of Lemma~\ref{Parab} for $g_{N }$ with the obvious bound $ \| e^{-tH}\|\leq e^{-\lambda_{1}t}$, one can estimate function \e{eq:par7f} for small $t$ by $C_{N }   t^{N-d/4+1}$ uniformly in $x, x^\prime$ from compact sets. This yields the following result.

   \begin{theorem}\label{Parab2}
Let the assumptions of   Lemma~\ref{Parab} hold.
Assume   that the operator $H=-\Delta +v(x)$ is self-adjoint on a domain ${\cal D}(H)\supset     {\cal S} ({\Bbb R}^d)$ and that it is semibounded from below. 
   Define   the functions $g_{n}$      by the recurrent relations $g_{0}(x,x^\prime)=1$ and \e{eq:par6}.   Then  
 \begin{equation}
G (x,x^\prime ;t)= (4\pi t)^{-d/2}   \exp (-  (4t)^{-1} |x - x^\prime |^2) \sum_{n=0}^\ii g_{n}(x,x^\prime)t^n 
\label{eq:Par}\end{equation}
as $t\ri 0$ uniformly with respect to $x$ and $x^\prime$ from   compact subsets of ${\Bbb R}^d$. Moreover, asymptotic expansion \e{eq:Par} can be infinitely differentiated in the variables
$x$, $x^\prime$ and $t$. In particular, we have that
 \begin{equation}
G (x,x  ;t)= (4\pi t)^{-d/2}     \sum_{n=0}^\infty g_{n}(x  )t^n, \q g_{n}(x)=g_{n}(x,x ),\q  t\ri 0.
\label{eq:Par5}\end{equation}
 \end{theorem}
 
 Of course expansion \e{eq:Par} is of interest for $|x - x^\prime |=O(t^{1/2})$ only.
 
  Proceeding from  \e{eq:par6}, one can   give  closed expressions  for the functions $g_{n} (x,x^\prime)$. In particular, it is not difficult to calculate explicitly  the first functions $g_{n}(x)$ (known as local heat invariants of the operator $H$):
  \[
     g_{1} =-v , \q g_{2} =2^{-1}v^2  -6^{-1} \Delta v , \q 
     g_{3} = -6^{-1}(v^3  - v \Delta v  -2^{-1} |\nabla v |^2+ 10^{-1}\Delta^2 v),
     \]
     \begin{eqnarray*} 
     g_{4} =24^{-1}v^4 + 30^{-1} \langle \nabla v, \nabla   (\Delta v) \rangle +
  60^{-1}   v \Delta^2 v +72^{-1}   (\Delta  v)^2-   840^{-1}  \Delta^3 v
  \\
  -   12^{-1} v^2 \Delta  v  -   12^{-1} v  |\nabla v |^2        +   90^{-1}   \tr (\hess v )^2.
   \end{eqnarray*} 
  For an arbitrary $n$,  the functions $g_{n}(x)$ can be found  (see \cite{Pol}  and also \cite{HiPol} where the results of \cite{AKa} were used) by the formula
  \begin{equation} 
g_{n}(x)= (-1)^n \Gamma(n+d/2)   
   \sum_{k=0}^{n-1}\frac{(-\Delta_{x }+ v(x))^{k+n} (|x-x^\prime|^{2k}) \Big|_{x =x^\prime} }{4^k k! (k+n)!(n-1-k)!\Gamma(k+d/2+1)}. 
  \label{eq:AsExpHe1}\end{equation}

  \medskip
  
  {\bf 3}.
Since
\begin{equation}
 R  (x, x^\prime; z) =\int_{0}^\infty  G(x, x^\prime; t)    e^{tz} dt, \q \Re z< \inf\sigma (H)=\lambda_{1},
\label{eq:loc}\end{equation}
we can relate 
 expansion \e{eq:Par} as $t\ri 0$ with the asymptotic expansion  of the resolvent kernel  $R(x ,x^\prime; z)$ as $|z|\ri \infty$ and thus to enhance the results of  subs.~1. We emphasize that   the boundedness of $v(x)$ is not required now.

 \begin{theorem}\label{Parab3}
Under   the assumptions of Theorem~$\ref{Parab2}$,   for sufficiently large $N$,
 \begin{equation}
R  (x,x^\prime ;z)=   \sum_{n=0}^N  g_{n}(x,x^\prime)  R_{0}^{( n)} (x, x^\prime; z)
+ O(|z|^{-N+d/2 - 2 }) 
\label{eq:loc1}\end{equation}
  as $|z|\ri \infty$, $\Re z < \inf \sigma(H)$, uniformly with respect to $x$ and $x^\prime$ from a compact subset of ${\Bbb R}^d$. Expansion \e{eq:loc1} can be infinitely differentiated in $x$, $x^\prime$ and $z$.
   \end{theorem}

Theorems~\ref{AsyRes} and \ref{Parab3} give two different expansions of the resolvent kernel.
Similarly to subs.~1, we can pass in \e{eq:loc1} to the limit $x^\prime\ri x$.
 Recall that kernel $R_{0}  (x, x^\prime; z)$ is   singular on the diagonal, but the derivatives $R_{0}^{( n)} (x, x^\prime; z)$ are getting smoother as $n$ increases. In particular, the function $R_{0}^{(n)}  (x, x^\prime ; z)$ is continuous on the diagonal if $2(n+1)>d$ and
\[
R_{0}^{(n)}  (x, x ; z)= (4\pi )^{-d/2} \Gamma(1+n -d/2) (-z)^{d/2- n -1}. 
\]
Therefore removing from \e{eq:loc1} $R_{0}  (x,x^\prime ;z)$ and differentiating $(m-1)$-times with respect to $z$, we can set $x= x^\prime$ if $2(m+1)>d$. This yields the AE
  \[
 (4\pi )^{-d/2}   \sum_{n=0}^\infty  \Gamma (n+m+1-d/2) g_{n+1} (x)  
  (-z)^{d/2 -m-n -1 }
  \]
  for the left-hand side of \e{eq:loc2a}.
Comparing it with the right-hand side of  \e{eq:loc2a}, we find that
      \begin{equation} 
 r_{n}^{(m)}  (x)=   (4 \pi )^{-d/2 }   \, (m-1)!^{-1}\Gamma ( n +m+1 -d/2) \, g_{n+1}  (x).
\label{eq:AsExpHegg}\end{equation}

Integrating \e{eq:Par5} over $x\in {\Bbb R}^d$, we formally obtain the AE
    \begin{equation} 
\tr ( e^{- H  t} -e^{- H_{0}  t})  =  (4 \pi t)^{-d/2 }\sum_{n=1}^\infty \mathbf{g}_{n}  t^n, \q
{\rm    where }\q  \mathbf{g}_{n}=\int_{{\Bbb R}^d} g_{n}(x) dx.
\label{eq:AsExpHe3}\end{equation}
The passage from  \e{eq:Par5}  to  \e{eq:AsExpHe3} requires   some estimates of
 $F_{N}(x,x ;t)$ at infinity which of course demand an appropriate decay of a potential as 
$|x|\ri\infty$. We avoid such estimates and  deduce AE  \e{eq:AsExpHe3}  from Theorem~\ref{AsExpTRMP} using  the formula (the  inversion of the Laplace transform \e{eq:loc})
\[
\tr(e^{- H  t}-e^{- H_{0}  t} )  = (2\pi i)^{-1} (m-1)!  t^{-m+1}  \int_{\beta-i\infty}^{\beta+i\infty}  
\tr(R^m(z)-R_{0}^m (z)) e^{-tz} d z
\]
where $2(m+1)>d$ and $\beta <\lambda_{1}$.   Let us formulate the precise result.
   
 \begin{theorem}\label{ReHe1}
 Suppose that  $v(x)$    satisfies estimates
\e{eq:1.4LR} where $\rho> d$.  
Then   AE    \e{eq:AsExpHe3}  holds.
   \end{theorem}

  We note finally that according to \e{eq:AsExpHegg}
   \begin{equation} 
 \mathbf{r}_{n}^{(m)}   =   (4 \pi )^{-d/2 }   \,(m-1)!^{-1} \Gamma ( n +m+1 -d/2) \, \mathbf{g}_{n+1} .
\label{eq:AsReS}\end{equation} 

     %%%%%%%%%%%%%%%%%%%%

\section{High-energy asymptotics of the SM}

%%%%%%%%%%%%%%%%%%%%%%%

{\bf 1.} 
Let assumption (\ref{eq:1.4LR}) be satisfied. Away from the diagonal $\omega = \omega^\prime$,
the integral kernel $s(\omega,\omega^\prime;\lambda)$ of the SM $S(\lambda)$ is  
$C^\infty$-function   and decays faster than any power of 
$\lambda^{-1}$ as $\lambda\ri\ii$.
On the contrary, it acquires diagonal singularities which are determined by the
fall-off of $v(x)$ at infinity. It turns out that  these singularities and   the
high-energy limit are described  by the same formulas.  

To describe $s(\omega,\omega^\prime;\lambda)$ in a neighbourhood of the diagonal, we recall first a standard construction of   approximate but explicit solutions $ \psi_N $
 of the  Schr\"odinger equation $-\Delta\psi +v(x)\psi=k^2 \psi $. This construction relies on a solution of the corresponding  transport equation  by iterations.
  Let us set
\[ 
 \psi_N (x,\xi)=e^{i\langle x,\xi\rangle} {\rm b}_N  (x,\xi),\quad
\xi=k \omega\in{\Bbb R}^d, \q k=\sqrt{\lambda},
\]
where
\begin{equation}
 {\rm b}_N(x,\xi)= \sum_{n=0}^N (2ik)^{-n} b_n(x, \omega),\quad b_0(x,
\omega)=1. 
\label{eq:WF1}\end{equation}
 Plugging these  expressions into the Schr\"odinger equation   and equating coefficients at the same powers of
$(2ik)^{-1} $, we obtain recurrent equations  
\[
   \langle \omega,\nabla_x b_{n+1}(x, \omega)\rangle = -\Delta b_n(x,\omega)  +v(x)
b_n(x,\omega)  
\]
whence
\begin{equation}
   b_{n+1}(x,\omega)  =\int_{-\infty}^0 \Bigl( -\Delta b_n(x+t \omega, 
\omega)+v(x+t \omega) b_n(x+t \omega, \omega)\Bigr)dt.
\label{eq:4BBB}\end{equation}
It is easy to see that under  assumption \e{eq:1.4LR} for any $c<1$ and $\varepsilon_{0}=\min\{1, \rho-1\}$
\begin{equation}
 |\partial_x^\alpha    b_n(x,\omega)| \leq C_{\alpha } 
(1+|x|)^{-\rho+1-\varepsilon_{0}(n-1) -|\alpha| },\q n\geq 1,   \q {\rm if} \q \langle x , \omega \rangle\leq   c | x|.
\label{eq:WF7}\end{equation}

 Let us fix some point
$\omega_0\in {\Bbb S}^{d-1}$. Let $\Lambda (\omega_0)$ be the plane orthogonal to
$\omega_0$ and $\Omega=\Omega  (\omega_0,\delta)\subset {\Bbb S}^{d-1}$ be determined by the condition
$  \langle\omega ,\omega_0\rangle>\delta>0$.  Set
$ x=\omega_0 z+y$, where $ y\in \Lambda  (\omega_0)$,
\begin{eqnarray} 
A_N (\omega,\omega^{\prime}, y;\lambda)= 2^{-1}  \langle\omega+\omega^\prime,\omega_0\rangle {\rm b}_{N} ( y,-k\omega) {\rm b}_{N} ( y,
k\omega^\prime)
\nonumber\\
 +(2i k)^{-1 }\Bigl( {\rm b}_{N} (y,-k\omega ) (\partial_z {\rm b}_{N} )( y,
k\omega^\prime)-{\rm b}_{N} ( y,k\omega^\prime) (\partial_z
{\rm b}_{N} )(  y,-k\omega)\Bigr) 
\label{eq:SST}\end{eqnarray}
 and
\begin{equation} 
 s_N (\omega,\omega^\prime ;\lambda)=
   (2\pi)^{-d+1} k^{d-1}
  \int_{\Lambda_{\omega_0}}e^{ik  \langle y ,\omega^\prime-\omega \rangle}
A_N(\omega,\omega^{\prime}, y;\lambda) dy 
\label{eq:SSTX}\end{equation}
 for $\omega,\omega^\prime\in \Omega $.  Since $| \langle y, \omega  \rangle |
\leq c |y |$ and   $| \langle y, \omega^\prime \rangle |
\leq c |y|$ where $c<1$ 
 for $\omega ,\omega^\prime \in \Omega $, $y\in \Lambda  (\omega_0)$,  estimates
(\ref{eq:WF7}) imply that oscillating 
 integral (\ref{eq:SSTX}) is well defined. As shown  in \cite{Yind},
 the function $ s_N (\omega,\omega^\prime ;\lambda)$ describes all singularities of 
$s(\omega,\omega^\prime;\lambda)$ and approximates it with arbitrary accuracy as
$\lambda\rightarrow\infty$.

\begin{theorem}\label{AsExp}
 Let assumption $(\ref{eq:1.4LR})$ hold for $\rho> 1$.   Let the function ${\rm b}_N(x,\xi)$ be defined by formulas $(\ref{eq:WF1})$ and
$(\ref{eq:4BBB})$. Define for $\omega,\omega^\prime\in\Omega $ the function
$ s_N$ by equalities  $(\ref{eq:SST})$, \e{eq:SSTX}. Then,  for     an arbitrary number $p$ and a sufficiently large $N=N(p)$, the remainder
$
 s (\omega,\omega^\prime;\lambda)-s_N(\omega,\omega^\prime ;\lambda) 
$
belongs to the class $C^p (\Omega\times\Omega)$ and the $C^p$-norm of this function is
$O(\lambda^{-p})$ as $\lambda\rightarrow\infty$.
\end{theorem}
  
  Another form of AE  of the scattering amplitude
  $ s (\omega,\omega^\prime,\lambda)$ can be found in \cite{Skr}. Theorem~\ref{AsExp}
  extends to long-range potentials \cite{Ylr}.

 Formula   (\ref{eq:SSTX}) shows  that we actually consider  the SM $S(\lambda)$ as a pseudo-differential operator (on
the unit sphere) determined by its amplitude $ A  (\omega,\omega^{\prime}, y;\lambda)$. Plugging expression (\ref{eq:WF1})  into  the right-hand side of (\ref{eq:SST}), we can expand this amplitude into asymptotic series
in powers of $(2 ik)^{-1}$. 
Moreover, it is possible to reformulate Theorem~\ref{AsExp}  using the standard procedure (see, e.g., \cite{Sh})  of passage from an amplitude  of a pseudo-differential operator to its  symbol. 
We note also that the operator-valued function $S\in C^\ii ({\Bbb R}_{+})$,  and its AE as $\lambda\ri \ii$ can be infinitely differentiated with respect to
 $\lambda$.

\medskip

{\bf 2.} 
Next  we consider the operator
$
T(\lambda)=-i S^*(\lambda) S^{\prime}(\lambda)  
$
 known as the Eisenbud-Wigner time-delay operator. In view of unitarity of the SM, the operator $T(\lambda)$ is self-adjoint. Remark that its kernel
$
t(\omega,\omega^\prime;\lambda) =O(\lambda^{-\infty}) \q\mathrm{for} \q\omega\neq\omega^\prime.
$
Standard results of the pseudo-differential calculus (see, e.g., \cite{Sh})  show that, together with $S(\lambda)$ and $ S^{\prime}(\lambda)$, the operator  $T(\lambda)$ is also a pseudo-differential operator on the unit sphere. Moreover, Theorem~\ref{AsExp} entails
  the following result which we formulate in terms of its right symbol.

  \begin{theorem}\label{AsExpT}
 Let assumption $(\ref{eq:1.4LR})$  hold for $\rho>1$. Then kernel of the time-delay operator $T(\lambda)$ admits expansion into the asymptotic series
 \begin{equation} 
 t ( \omega, \omega^{\prime};\lambda) =     (2\pi)^{-d+1} k^{d-4} 
  \sum_{n=0}^\infty (2ik)^{-n}
\int_{\Lambda_{\omega }} e^{ik  \langle y ,\omega^{\prime} -\omega   \rangle} {\tt t}_{n}   (  \omega^{\prime},  y ) dy
\label{eq:SSTX3t}\end{equation} 
where ${\tt t}_n  = \bar{\tt t}_{n}$ are smooth functions of  $  \omega^{\prime}\in \Omega (\omega ,\delta)$, $y\in  \Lambda (\omega)$ and $\partial_{y}^\alpha {\tt t}_{n} (  \omega^{\prime} , y)= O(|y|^{-\rho+1-\varepsilon_{0}n-|\alpha|})$ as $|y|\ri\ii$.
If $\rho>d$, then $t$ is a continuous function of $\omega, \omega^{\prime}$ and
  \begin{equation} 
 t ( \omega, \omega;\lambda) =  (2\pi)^{-d+1} \lambda^{d/2-2} \sum_{n=0}^\infty 
  (-   4  \lambda)^{-n}
\tau_{n}  ( \omega ),
\label{eq:SSTX3td}\end{equation} 
where   $  \tau_{n}   ( \omega )$ are the integrals of $ {\tt t}_{2n}  ( \omega ,  y )$ over $y\in \Lambda (\omega)$.  Moreover, $T \in C^\ii ({\Bbb R}_{+})$,  and
  AE \e {eq:SSTX3t}, \e {eq:SSTX3td} can be infinitely differentiated with respect to $\lambda$.
 \end{theorem}

%%%%%%%%%%%%%%%%%%%%%%%%%%%%%%%%%%%
%%%%%%%%%%%%%%%%%%%%%%%%%%%%%%%%%%%
\section{High-energy asymptotics of the SSF}
 %%%%%%%%%%%%%%%%%%%%%%%%%%%%%%%%%%%%%
%%%%%%%%%%%%%%%%%%%%%%%%%%%%%%%%%%%

Our goal in this section is to find the complete AE      of the SSF $\xi (\lambda)$ as $\lambda\ri\infty$. At the same time we extend AE (\ref{eq:AsExpTRK7})   to the whole complex plane cut along the positive half-axis. In subs.~1 we prove the existence of the complete AE  of  
 $\xi (\lambda)$. However the method of this subsection   gives   complicated expressions for coefficients $\xi_{n}$ of this expansion. This drawback is remedied in subs.~2 and 3 where a connection of the results of subs.~1 with the results of   subs.~1 and 2, \S 4, is established. This yields expressions for $\xi_{n}$ in terms of the coefficients 
${\bf r}_{n}^{(m)}$ as well as of the heat invariants ${\bf g}_{n+1} $.

\medskip

{\bf 1.}
We proceed from equation (\ref{eq:BKB}) which follows from   the Birman-Kre\u{\i}n formula  (\ref{eq:BK}) because the SSF $\xi$ is a continuous function of $\lambda>0$   and the operator-valued function $S (\lambda)$ is differentiable with respect to $\lambda$ in the trace norm.
  Let   us now use   Theorem~\ref{AsExpT}.
Integrating \e{eq:SSTX3td} over $\omega\in {\Bbb S}^{d-1}$ and taking into account formula \e{eq:BKB}, we find the AE of the derivative $\xi^\prime (\lambda)$ of the SSF.

 \begin{theorem}\label{AsExpSSF}
 Let assumption $(\ref{eq:1.4LR})$ where $\rho>d$ hold. Then $\xi \in C^\ii ({\Bbb R}_{+})$ and the   AE
  \begin{equation} 
\xi^\prime (\lambda)=    \lambda^{d/2-2} \sum_{n=0}^\infty  \eta_{n} \lambda ^{-n}
\label{eq:SSFX3t}\end{equation} 
holds with the coefficients
\begin{equation}
 \eta_0      =  4^{-1} (d-2) (2\pi)^{-d }  |{\Bbb S}^{d-1}|\int_{{\Bbb R}^{d }} v( x) dx,
 \q \eta_{n}  =- (2\pi)^{-d } (-4 )^{-n}\int_{{\Bbb S}^{d-1}}\tau_{n}  ( \omega )d\omega.
\label{eq:SSFX2t}\end{equation}
 Expansion \e {eq:SSFX3t}  can be infinitely differentiated.
\end{theorem}

Integrating \e{eq:SSFX3t} over $\lambda$, we obtain also the following result.

 \begin{corollary}\label{AsExpSSFc}
For a suitable constant $\gamma$,  the  SSF $\xi(\lambda)$ admits expansion into the asymptotic series
  \begin{equation} 
\xi  (\lambda)=   \lambda^{d/2-1} \sum_{n=0, n\neq d/2-1}^\infty
   \xi_{n}  \lambda^{-n}
 +   \tilde{\xi} \ln\lambda +\gamma,
\label{eq:SSFX3tc}\end{equation}
 where  $\xi_{n}= (d/2-n-1)^{-1} \eta_{n}$,   $\tilde{\xi}= 0$ for $d$ odd and   $\tilde{\xi}=\eta_{d/2 -1}$ for $d$ even.
\end{corollary}

Expressions \e{eq:SSFX2t} for the coefficients $\eta_{n}$, $n\geq 1$, are rather complicated. Below we shall find efficient expressions   for   $\xi_{n}$.
  At the same time we shall show that $\tilde{\xi}=0$ for  all $d$ and  find the constant $\gamma$ in  \e{eq:SSFX3tc}.

\medskip

 {\bf 2.}
 In view of representation \e{eq:ssfd0}, the following assertion can be
 deduced from Corollary~\ref{AsExpSSFc} with a help of a version  of  the Privalov theorem.   
 
 \begin{proposition}\label{AsExpTR}
 Let assumption $(\ref{eq:1.4LR})$ where $\rho>d$ hold, and let $2(m+1)>d$. Then,
 for some numbers $c^{(m)}_{n}$, $n\geq m+1$, the expansions into asymptotic series hold as $|z|\ri\infty$, $z\in{\Bbb C}\setminus {\Bbb R}_{+}:$
 \begin{equation} 
\tr(R^m(z)-R^m_{0}(z))=   \sum_{n=0}^\infty {\bf r}_{n}^{(m)} (-z)^{d/2-m-n-1 } -\gamma (-z)^{-m} +\sum_{n=m+1}^\infty c_{n}^{(m)}z^{-n}
\label{eq:AsExpTR2}\end{equation}
   if $d$ is odd, and
  \begin{eqnarray} 
\tr(R^m(z)-R_{0}^m(z)) =    \sum_{n=0}^{d/2-2} {\bf r}_{n}^{(m)} (-z)^{d/2-m-n-1 }
  \nonumber\\
+   \sum_{n=d/2 }^\infty \sigma_{n}^{(m)} (-z)^{d/2-m-n-1 }\ln(-z)  -\gamma (-z)^{-m}  
   +\sum_{n=m+1}^\infty c_{n}^{(m)} z^{-n}  ,
\label{eq:AsExpTR1}\end{eqnarray} 
  if $d$ is even $($the sum over $n=0,\ldots, d/2-2$ is absent if $d=2)$.
  Here $\gamma$ is the same as in \e{eq:SSFX3tc} and the asymptotic coefficients 
  $ {\bf r}_{n}^{(m)} $ are defined by the formula
     \begin{equation} 
  {\bf r}_{n}^{(m)} =- m \mathsf{B}(-n +d/2, m+n+1-d/2)\,  \xi_{n}
\label{eq:AsExpTRcoef}\end{equation}
 $( \mathsf{B}(p,q)$   is the beta-function$) $
for all $n$     if $d$ is odd and for $n=0,   \ldots, d/2-2$ if $d$ is even and
      $  \sigma_{n}^{(m)} = - (-1)^{n-d/2 }
  m \begin{pmatrix}
  n-d/2+m\\ m
   \end{pmatrix}
     \xi_{n}$ for $ n= d/2 , d/2+1 ,\ldots$.
\end{proposition}

We emphasize that the coefficients $c_{n}^{(m)}$ are not determined by the AE of the SSF $\xi(\lambda)$ as $\lambda\ri\ii$.
 Let us now compare Proposition~\ref{AsExpTR} with Theorem~\ref{AsExpTRMP}.
           Expansion (\ref{eq:AsExpTRK7})  should coincide with
\e{eq:AsExpTR2} or \e{eq:AsExpTR1} in the  angle  $z\in\Pi_{\theta}$. Therefore 
the coefficients  (\ref{eq:AsExpTRK6}) and \e{eq:AsExpTRcoef} are the same.
Moreover, $\gamma=0$, $c_{n}^{(m)} =0$ in 
\e{eq:AsExpTR2} and $\gamma=- {\bf r}_{d/2-1}^{(m)}$, $\sigma_{n}^{(m)}=0$,  $c_{n}^{(m)}={\bf r}_{n+d/2-1}^{(m)}$ in \e{eq:AsExpTR1}. 
Thus,   Proposition~\ref{AsExpTR} implies the following result.

  \begin{theorem}\label{AsExpTRM}
    Let assumption $(\ref{eq:1.4LR})$ with $\rho>d$  hold, and let $2(m+1)>d$.
      Then    AE $(\ref{eq:AsExpTRK7})$ 
      as $|z|\ri\infty$   is valid for all  $z\in {\Bbb C}\setminus {\Bbb R}_{+}$. 
\end{theorem}

At the same time we have improved significantly  Theorem~\ref{AsExpSSF}  and Corollary~\ref{AsExpSSFc} obtaining explicit expressions of coefficients $\xi_{n} $ in AE \e{eq:SSFX3tc}  in terms of ${\bf r}_{n}^{(m)}$.

 \begin{theorem}\label{AsExpSSF2}
    Let assumption $(\ref{eq:1.4LR})$ with $\rho>d$   hold.
 Then   the  SSF $\xi(\lambda)$ admits expansion into the asymptotic series
  \begin{equation} 
\xi  (\lambda)=   \lambda^{d/2-1} \sum_{n=0 }^\infty
    \xi_{n}  \lambda^{-n},
\label{eq:SSFodd}\end{equation}
    where   the coefficients $\xi_{n} $ are determined by formulas \e{eq:AsExpTRcoef} and $(\ref{eq:AsExpTloc})$, $(\ref{eq:AsExpTRK6})$. 
        In particular, $\xi_{n}=0$ for all $n\geq d/2$ if $d$ is even.
         \end{theorem}

 \medskip

{\bf 3.} 
Alternatively, AE  \e{eq:SSFodd} can be deduced from Theorem~\ref{ReHe1} (instead of Theorem~\ref{AsExpTRMP}).    Since
  \[ 
\tr ( e^{- H  t} -e^{- H_{0}  t})= -t \int_{-\infty}^\infty\xi  (\lambda) e^{-t\lambda} d\lambda,
\]
   it   follows (see, e.g., \cite{Col})  from AE  \e{eq:SSFX3tc} of the SSF that
  \begin{equation} 
\tr ( e^{- H  t} -e^{- H_{0}  t}) = -\sum_{n=0}^\infty    \Gamma (d/2-n ) \xi_{n}  t^{n+1-d/2}    -\pi \gamma  +   \sum_{n=1}^\infty (-1)^n c_{n}   t^{n}  
   \label{eq:SsfH2}\end{equation}
   if $d $ is odd, and
 \begin{eqnarray} 
\tr ( e^{- H  t} -e^{- H_{0}  t})= -\sum_{n=0}^{d/2-1}    \Gamma (d/2-n )  \xi_{n} t^{n+1-d/2}  -\pi \gamma
 \nonumber\\
 -\sum_{n=d/2}^{ \infty }   (-1)^{n-d/2}    ( n-d/2 )!^{-1} \xi_{n}  t^{n+1-d/2} \ln t +   \sum_{n=1}^\infty (-1)^n c_{n}   t^{n} 
    \label{eq:SsfH3}\end{eqnarray}
   if $d $ is even. Here
     \begin{equation} 
c_{n}=   (n-1) !^{-1} \Big(\int_{-\infty}^0 \xi  (\lambda) \lambda^{n-1} d\lambda
    + \int_0^{\infty}  (\xi  (\lambda) -\sum_{ j < n+d/2-1} \xi_{j }\lambda^{d/2-j -1})\lambda^{n-1} d\lambda\Big).
  \label{eq:SHInt}\end{equation}
   
   Comparing \e{eq:SsfH2} and \e{eq:SsfH3} with (\ref{eq:AsExpHe3}), we recover Theorem~\ref{AsExpSSF2}. Moreover, we obtain explicit expressions
      \begin{equation} 
  \xi_{n} =- (4\pi)^{-d/2} \Gamma^{-1}(d/2 -n  ) \mathbf{g}_{n+1}, \q n=0,1, \ldots,  
\label{eq:xigg}\end{equation}
 for the asymptotic coefficients in \e{eq:SSFodd}
    in terms of the heat invariants  defined by formulas  \e{eq:AsExpHe1} and \e{eq:AsExpHe3}. Relation \e{eq:xigg} can also be obtained by putting together formulas (\ref{eq:AsReS}) and \e{eq:AsExpTRcoef}.

%%%%%%%%%%%%%%%%%%%%%%%%%%%%%%%%%%%
%%%%%%%%%%%%%%%%%%%%%%%%%%%%%%%%%%%
\section  {Trace identities  }
 %%%%%%%%%%%%%%%%%%%%%%%%%%%%%%%%%%%%%
%%%%%%%%%%%%%%%%%%%%%%%%%%%%%%%%%%%

The contribution of the continuous spectrum to 
trace identities  is, roughly speaking, given in terms of the function $\tr (R(\lambda+i0) - R_{0}(\lambda+i0))$. In the case of integer order the imaginary part and in the case of half-integer order the real part of this function are considered.

\medskip

{\bf 1.}
Let us start with formulas of integer order.  The method used here does not require any specific study of the PD or   the SSF in a neighborhood of the point $\lambda=0$. In particular, in the case $d=1$ Theorem~\ref{AsExpSSFTr1} does not exclude that the operator $H$ has   infinite number of negative eigenvalues.

  Let us   compare the coefficient  at $t^n$, $n=1,2,\ldots$, in expansion (\ref{eq:AsExpHe3})   with the same coefficient   in expansions \e{eq:SsfH2} or \e{eq:SsfH3}. According to (\ref{eq:AsExpHe3})   this coefficient is zero if $d$ is odd and it equals 
  $(4\pi )^{-d/2}\mathbf{g}_{d/2+n }$  if $d$ is even.   On the other hand, according to \e{eq:SsfH2} or \e{eq:SsfH3} it equals $(-1)^n c_{n}$ and hence   is   determined by  equality  \e{eq:SHInt}. Let us also take into account that, for $\lambda<0$, the SSF is given  by formula \e{eq:5.3.7D2x} which yields
  \[
  \int_{-\infty}^0   \xi(\lambda) \lambda^{n-1}d\lambda= n^{-1} \sum_{j=1}^{\mathsf{N}} \lambda_{j}^{n}.
  \]
  Thus, the equalities $c_{n } = 0$ for $d$ odd and $c_{n} = (-1)^n  (4\pi )^{-d/2}\mathbf{g}_{d/2+n }$ for $d$ even give us identities known as the trace identities.
  
\begin{theorem}\label{AsExpSSFTr1}
    Let assumption $(\ref{eq:1.4LR})$ where $\rho>d$ hold.
  Then   for all $n=1,2,\ldots$, we have the following identities. If $d$ is odd, then
        \[
 \int_{0}^\infty  \Big(\xi(\lambda)- \sum_{j=0}^{(d-3)/2+n} \xi_{j} \lambda^{d/2-j-1}  \Big)\lambda^{n-1}d\lambda+ n^{-1}  \sum_{j=1}^{\mathsf{N}} \lambda_{j}^{n}=0.
\]
 If $d$ is even, then
       \[ 
       \int_{0}^\infty  \Big(\xi(\lambda)- \sum_{j=0}^{d/2-1 } \xi_{j} \lambda^{d/2-j-1}  \Big)\lambda^{n-1}d\lambda+ n^{-1}  \sum_{j=1}^{\mathsf{N}} \lambda_{j}^{n}
       = (-1)^n(4\pi )^{-d/2}(n-1)! \mathbf{g}_{d/2+n }. 
       \]
    \end{theorem}

\medskip

{\bf 2.}
Next we consider trace identities of half-integer order. Now we suppose that $d\leq 3$ which allows us to conveniently formulate results in terms of the regularized PD $D_{2}(z)$     discussed in  \S 3 (for $d=1$  the results can be formulated in terms of $D (z)$ - see \cite{FaZa}).  
%Recall that $D_{2}(z)$ is an analytic function of $z$, and 
%the branch of $\ln D_{2}(z)$ is fixed by condition \e{eq:5.arg}. 

 \begin{proposition}\label{AsExpSSFTrh}
    Let $d\leq 3$, and let assumption $(\ref{eq:1.4LR})$  where $\rho> d/2$ hold.
      Then  the AE  
              \begin{equation} 
\ln D_{2} (z)  = \sum_{n=1}^\infty   \delta_{n}
  (-z)^{d/2-1-n  }     
\label{eq:AsExpSSFTrh1}\end{equation}
              as $|z|\ri\infty$, $z\in {\Bbb C}\setminus {\Bbb R}_{+}$,
 is valid with the coefficients
$  \delta_{n} = - (n+1-d/2)^{-1}  {\bf r}_{n}^{(1)}$
where the coefficients ${\bf r}_{n}^{(1)}$ are  defined by 
\e{eq:AsExpTloc} and \e{eq:AsExpTRK6}.
  In particular, as $\lambda\ri\infty$, 
$\ln |D_{2} (\lambda\pm i0) | = O(\lambda^{-\infty})$ 
  if $d$ is odd, and
  \[ 
\ln |D_{2} (\lambda\pm i0) | =  (-1)^{d/2 +1} \sum_{n=1}^\infty  (-1)^n \delta_{n}
  \lambda^{d/2-1-n  }     
\]
  if $d$ is even.
\end{proposition}

Indeed, let us  use Remark~\ref{AsExpTRMPr} to Theorem~\ref{AsExpTRMP}. Combining equation (\ref{eq:PDRp}) for $p=2$  with AE  (\ref{eq:AsExpTRK7r}), we see that
\[
d \ln D_{2} (z)/dz  = - \sum_{n=1}^\infty   {\bf r}^{(1)}_{n}
  (-z)^{d/2-2-n  }.    
\]
  It remains to integrate it and to  take   relation \e{eq:5.3.2}  into account.

By virtue of (\ref{eq:AsReS}), we also have   that
 $  \delta_{n} = - (4\pi )^{-d/2} \Gamma (n+1-d/2)  {\bf g}_{n+1}$.

\begin{theorem}\label{AExpSSFTrh1} 
 Let $d\leq 3$,   let assumption $(\ref{eq:1.4LR})$ where $\rho>2$  be satisfied
 and let $n=0,1,2,\ldots$.
 Suppose also that condition \e{eq:7.3.15} is satisfied for $\alpha=n+1/2$.
  Let the numbers $\delta_{n}$ be the same as in Proposition~$\ref{AsExpSSFTrh}$.
  Then, for $d=1, 3$,    
          \[ 
\pi^{-1} (-1)^n \int_{0}^\infty  \ln | D_{2}(\lambda+i 0) | \lambda^{n-1/2}d\lambda - (n +1/2)^{-1} \sum_{j=1}^{\mathsf{N}} |\lambda_{j}|^{n+1/2}=  \delta_{ n +1}
\]
        and,  for $d=2$  $($the sum under the integral sign is absent for $n=0)$,
       \begin{eqnarray*} 
\pi^{-1} (-1)^n  \int_{0}^\infty  \Big(\ln |D_{2}(\lambda+i 0)| - \sum_{j=1}^{n} (-1)^{j}\delta_{j} \lambda^{d/2-1-j}  \Big)\lambda^{n-1/2}d\lambda
       \nonumber\\
        - (n+1/2)^{-1} \sum_{j=1}^{\mathsf{N}} |\lambda_{j}|^{n+1/2}=0.
        \end{eqnarray*}
      \end{theorem}
    
    The proof is similar to that of   Theorem~\ref{5.3.2}.
    By virtue of the Cauchy theorem,   we have that
  \begin{equation} 
 \int_{\Gamma_{R,\varepsilon}} \ln D_{2} (z)( -z)^{n-1/2} dz=0  
\label{eq:AExpSSFTrh2}\end{equation}
where $\Gamma_{R,\varepsilon}$ is the same contour as in \e{eq:5.2.23}.  The integrals
 over $C^\pm_\varepsilon(\lambda_j)$ and $C^\pm_\varepsilon(0)$ tend again  to zero as $\varepsilon\rightarrow 0$.  Thus, it follows from \e{eq:AExpSSFTrh2}  that
 \begin{eqnarray}  
2  (-1)^n  \int_{0}^R \ln | D_{2} (\lambda + i0)|
   \lambda^{n-1/2} d\lambda -2\pi   (n+1/2)^{-1} \sum_{j=1}^\mathsf{N} |\lambda_{j}|^{n+1/2}
   \nonumber\\
   = i\int_{C_{R }} \ln D (z)( -z)^{n-1/2} dz . 
\label{eq:AxpSSFTrh7}\end{eqnarray}
The AE  as $R\ri\infty$ of both integrals  
  can   easily be deduced from \e{eq:AsExpSSFTrh1}. Hereby the divergent terms cancel each other.
It remains to pass   to the limit $R\ri\infty$ in \e{eq:AxpSSFTrh7}.

 Trace identities of integer order can also be obtained by the method of proof of Theorem~\ref{AExpSSFTrh1} (hereby $n-1/2$ should be replaced by $n-1$ in \e{eq:AExpSSFTrh2}) which however requires more stringent assumptions than Theorem~\ref{AsExpSSFTr1}.

\end{document}